\documentstyle[12pt]{article}
 \newcommand{\be}{\begin{equation}}
 \newcommand{\ee}{\end{equation}}
 
\newcommand{\bC}{{\bf C}}

 \newcommand{\cA}{{\cal A}}
 
 \newcommand{\cC}{{\cal C}}

 \newcommand{\cF}{{\cal F}}

 \newcommand{\cO}{{\cal O}}

 \newcommand{\ra}{\rightarrow}
 \newcommand{\lra}{\longrightarrow}
 
 \newcommand{\kahler}{K\"{a}hler}

 \newcommand{\im}{\mathop{{\fam0 Im}}\nolimits}

 \newcommand{\Ker}{\mathop{{\fam0 Ker}}\nolimits}
 
 \newcommand{\Rep}{\mathop{{\fam0 Rep}}\nolimits}
 \newcommand{\Res}{\mathop{{\fam0 Res}}\nolimits}
 \newcommand{\aut}{\mathop{{\fam0 Aut}}\nolimits}
 \newcommand{\Int}{\mathop{{\fam0 Int}}\nolimits}
\newcommand{\ind}{\mathop{{\fam0 ind}}\nolimits}
 \newcommand{\End}{\mathop{{\fam0 End}}\nolimits}
  
 \newcommand{\rank}{\mathop{{\fam0 rank}}\nolimits}

 \newcommand{\tr}{\mathop{{\fam0 Tr}}\nolimits}

 \newcommand{\GL}{\mathop{{\fam0 GL}}\nolimits}
 
 \newcommand{\U}{\mathop{{\fam0 U}}\nolimits}

\newcommand{\C}{\mathop{{\fam0 c}}\nolimits}
\newcommand{\ch}{\mathop{{\fam0 ch}}\nolimits}

 \newcommand{\Vect}{\mathop{{\fam0 {Vect}}}\nolimits}

 \newcommand{\dbar}{\overline{\partial}}
 \newtheorem{definition}{Definition}[section]
 \newtheorem{lemma}[definition]{Lemma}
 \newtheorem{prop}[definition]{Proposition}
 \newtheorem{thm}[definition]{Theorem}
 \newtheorem{cor}[definition]{Corollary}

 \newcommand{\pf}{{\em Proof}. }
 \newcommand{\remark}{\noindent{\em Remark}. }

 \newcommand{\twist}{{\ast}}

 \newcommand{\qed}{\hfill$\Box$}

\begin{document}
\title{Twisted Higgs bundles and the fundamental group\\ of compact
\kahler\ manifolds}
\author{O. Garc\'{\i}a--Prada \and S. Ramanan}
\date{}
\maketitle
%\date
%{\footnotesize\tableofcontents\frenchspacing\smallbreak}

\begin{abstract}

We study polystable Higgs bundles twisted by a line bundle over a
compact K\"ahler
manifold.  These form a Tannakian category when  the first and second
Chern  classes of the bundle are zero. In this paper we identify the 
corresponding Tannaka group 
in the case in which the line bundle is of finite
order. This  group is described in terms of the pro-reductive completion of
the fundamental group of the manifold, and the character associated
to the line bundle.


\end{abstract}

%%%%%%%%%%%%%%%%%%%%%%%%%%%%%%%%%%%%%
\section*{Introduction}\label{intro}
%%%%%%%%%%%%%%%%%%%%%%%%%%%%%%%%%%%%%

Let $(X,\omega)$ be a compact \kahler\ manifold. A {\em Higgs bundle} 
over $X$ 
consists of a holomorphic vector bundle $E\ra X$ together with a sheaf
homomorphism
$\varphi:E\ra E\otimes \Omega^1$, satisfying $\varphi\wedge\varphi=0$,
where $\Omega^1$ is the holomorphic cotangent bundle of $X$.
A notion of stability for Higgs bundles,
similar to that of vector bundles, was first introduced for a Riemann
surface by Hitchin in
\cite{H1}, where he  proved that this is 
a necessary and sufficient condition
for the existence of an irreducible solution to 
the so-called {\em self-duality equations} for for a Hermitian metric
on  $E$. His results were later extended to the general case
by  Simpson \cite{Si1}.
Of special interest is  the moduli space of stable Higgs bundles $E$ for
which
$$
\C_1(E).[\omega]^{\dim X - 1}=0\;\;\; \mbox{and} \;\;\;
\ch_2(E).[\omega]^{\dim X - 2}=0,
$$
where $\C_1(E)$ and $\C_2(E)$ are the first and second Chern classes
of $E$, and $\ch_2(E)=1/2\C_1(E)^2-\C_2(E)$.
 The  moduli space of such Higgs bundles
can be identified with the moduli space of irreducible flat complex
connections,  which in
turn
is in correspondence with the moduli of  complex irreducible
representations of the fundamental group of $X$.
This is proved by  using, on one hand, Hitchin's and Simpson's  existence theorem 
and,
on the other, a theorem of 
Donaldson \cite{D2} and Corlette \cite{C} 
on the existence of {\em harmonic} metrics on flat bundles.

In this paper we shall deal with a {\em twisted} version of Higgs bundles.
Let $L$ be a holomorphic line bundle over $X$. An $L$-{\em twisted 
Higgs bundle} over $X$ (or just
twisted if there is no confusion), 
is a pair 
consisting of a holomorphic vector bundle $E$ and a sheaf morphism 
$\theta: E\ra E\otimes L\otimes \Omega^1$, i.e.
an element $\theta\in H^0(\End E\otimes L\otimes \Omega^1)$,
satisfying $\theta\wedge \theta=0$.
Stability is defined exactly in the same way as for  ordinary Higgs
bundles. 
There are many 
interesting aspects of Higgs bundles that are shared by the twisted 
 theory: twisted Higgs bundles, like ordinary Higgs bundles, have  nice
moduli spaces \cite{N,Y}, 
and define also, under certain conditions, complete integrable systems in
a generalised sense 
\cite{H2,Bo,Ma}.
There is one aspect, however, that to our knowledge has not yet been
studied, namely, their relation to the fundamental group.
The main indication that a certain relationship has to exist comes from
the fact that polystable $L$-twisted Higgs bundles, with the vanishing
condition on the first and second Chern characters as above,
like Higgs bundles, 
form a 
Tannakian category \cite{Si2}, which, by the Tannaka
duality theorem, is dual to a certain pro-reductive
group, from the representations of which the category can  be recovered.
In the ordinary Higgs bundle theory the group is 
the  pro-reductive completion of the  fundamental group of $X$
--- as one deduces from
the theorems of Hitchin, Simpson,  Donaldson, Corlette mentioned above.
The problem we wish to address is to find the corresponding  group in
the twisted situation.  In this paper, we take the first steps in this 
direction by considering the case in which the $\deg L=0$ and  $L$  is of
finite order, i.e. $L^n$ is isomorphic to the trivial line bundle for
some $n$.
In this situation $L$ corresponds to a unitary character 
$\pi_1(X)\ra \U(1)$ of the fundamental group of $X$, and we can give a 
description  of the group
in terms of the pro-reductive completion of $\pi_1(X)$ --- the Tannaka
group for the untwisted category --- and this character.

%%%%%%%%%%%%%%%%%%%%%%%%%%%%%%%%%%%%%%%%%%%%%%%%%%%%%%%%%%
\section{Twisted Higgs bundles}\label{twist}
%%%%%%%%%%%%%%%%%%%%%%%%%%%%%%%%%%%%%%%%%%%%%%%%%%%%%%%%%%
%%%%%%%%%%%%%%%%%%%%%%%%%%%%%%%%%%%%%%%%%%%%%%%%%%%%%%
\subsection{Twisted Higgs bundles and Hermitian metrics}
%%%%%%%%%%%%%%%%%%%%%%%%%%%%%%%%%%%%%%%%%%%%%%%%%%%%%%

Let $(X,\omega)$ be a compact \kahler\ manifold
and let $L$ be a holomorphic 
line bundle over $X$. An $L$-{\em twisted Higgs bundle} over $X$
is a pair $(E,\theta)$ consisting of a holomorphic vector bundle $E$
over $X$ and a Higgs field $\theta\in H^0(\End E\otimes L\otimes
\Omega^1)$,  satisfying $\theta\wedge\theta=0$,
where $\Omega^1$ is the holomorphic cotangent bundle of $X$.
A twisted  Higgs bundle $(E,\theta)$ is said to be {\em stable} if and only if
$\mu(E')<\mu(E)$ for every proper coherent subsheaf $E'\subset E$ invariant under
$\theta$, i.e. $\theta(E')\subset E'\otimes \Omega^1$. Recall that the
slope of a $E'$ is defined as $\mu(E')=\deg(E')/\rank (E')$, where 
$\deg(E')=c_1(E').[\omega]^{\dim X - 1}$.

The notion of stability is related to the existence of a special Hermitian  metric
on $E$. More precisely:

\begin{thm}\cite{Li}\label{hk}
Let $(E,\theta)$ be an $L$-twisted Higgs bundle. Let us fix a
Hermitian metric on $L$. The existence of a
Hermitian  metric $h$ on $E$ satisfying 
\be
\Lambda F_h +\Lambda[\theta,\theta^\ast]=\lambda I,\label{twisted-he}
\ee
is equivalent to the polystability of $(E,\theta)$.
\end{thm}
Here $F_h$ is the curvature of the unique connection
compatible with the Hermitian metric as well as  the holomorphic structure 
on $E$, and $\Lambda$ is the contraction with the \kahler\ form.
The constant $\lambda$ is determined by the slope of $E$, and $I$ is
the identity endomorphism of $E$.
By $\theta^\ast$ we denote the adjoint of $\theta$ with 
respect to $h$ and the metric of $L$, and 
$[\theta,\theta^\ast]= \theta\theta^\ast + \theta^\ast\theta$ is the usual 
extension of the Lie bracket to forms with values in the algebra of
endomorphisms. By {\em polystability} we mean that $(E,\theta)$
is a direct sum of twisted Higgs bundles of the same slope
as $E$ (zero in this case). 

When $L$ is the trivial line bundle one has the ordinary Higgs bundles
theory studied by Hitchin \cite{H1} on Riemann surfaces and Simpson
\cite{Si1} in the higher dimensional case.
When $\theta=0$,\  (\ref{twisted-he}) reduces
to the Hermitian--Einstein equation and one obtains the theorem 
of Narasimhan and Seshadri, Donaldson, and  Uhlenbeck and Yau 
\cite{AB,D1,D3,D4,NS,UY}.

Equation (\ref{twisted-he}) has, as in the untwisted case, a symplectic 
interpretation. It corresponds to the moment map for the action of
the unitary group on the product \kahler\ manifold
$\cA\times\Omega^{1,0}(End E\otimes L)$, where $\cA$ is the space of 
unitary connections on the $C^\infty$ Hermitian vector bundle $(E,h)$.
The moduli space of stable $L$-twisted Higgs bundles is then 
obtained as a \kahler\ quotient inheriting in this way a \kahler\
structure. A construction   of the moduli of $L$-twisted Higgs bundles 
using Geometric Invariant Theory has been given by Nitsure \cite{N}
for Riemann surfaces and  Yokogawa \cite{Y} in higher dimensions.

%%%%%%%%%%%%%%%%%%%%%%%%%%%%%%%%%%%%%%%%%%%%%%%%%%%%%%%%%%%%%%%%%%%%%%%%
\subsection{Twisted Higgs bundles and Tannakian categories}\label{tannaka}
%%%%%%%%%%%%%%%%%%%%%%%%%%%%%%%%%%%%%%%%%%%%%%%%%%%%%%%%%%%%%%%%%%%%%%%%
%%%%%%%%%%%%%%%%%%%%%%%%%%%%%%%%%%%%%%%%%%%%%%%%%%%%%%%%%%%%%%%%%%%%%%%%

Twisted Higgs bundles can be regarded from the point of view of Tannakian
categories.  This is the  point of view taken by Simpson  \cite{Si2}
in his study of ordinary Higgs bundles,
and we will follow  his approach. (See also \cite{De,DMOS,Sa,T} for more 
details about Tannakian categories).
A {\em tensor category} is a category $\cC$ with a functorial binary 
operation
$\otimes:\cC\times\cC\ra\cC$. An {\em associative} and {\em commutative} 
tensor
category is a tensor category provided with additional natural 
isomorphisms expressing associativity and commutativity of the 
tensor product  that have to satisfy certain {\em canonical axioms}.
A {\em unit} $1$ is an object $1$ provided with natural isomorphisms 
$1\otimes V\cong V$ satisfying canonical axioms. 
A functor $\cF$ between associative and commutative categories with 
unit is a functor  provided with natural isomorphisms 
$\cF(U\otimes V)\cong \cF(U)\otimes\cF(V)$.
A {\em neutral Tannakian category} $\cC$ is an associative and commutative
tensor category with unit, which is {\em abelian}, {\em rigid} 
(duals exist),  $\End(1)=\bC$, and which is provided with an 
{\em exact, faithful fibre functor} $\cF: \cC\ra \Vect$, where $\Vect$ 
is the tensor category of complex, finite dimensional vector spaces.

If  $G$ is an affine group scheme over $\bC$
the category  $\Rep(G)$  of  complex representations of $G$ is a 
neutral Tannakian category. The fibre functor $\cF_G$  is given by
by sending a representation of $G$ to the underlying vector space.
The group $G$ is recovered as the group $G={\aut}^\otimes (\cF_G)$ of 
tensor automorphisms of the fibre functor.
The converse is given by the fundamental duality theorem of
Tannaka--Grothendieck--Saavedra (\cite{De, DMOS,Sa}).

\begin{thm}
Let $(\cC,\cF)$ be a neutral Tannakian category and let 
$G=\aut^\otimes (\cF)$ be the group of tensor automorphisms of the
fibre funtor. Then $(\cC,\cF)\cong (\Rep(G),\cF_G)$.
\end{thm}

We shall briefly describe the group $\aut^\otimes (\cF)$
to which we shall refer sometimes as the {\em Tannaka group}
of the Tannakian category $(\cC,\cF)$ (see  \cite{Si2} and 
the references mentioned above for a detailed account.) 
Let $\End(\cF)$ be the algebra of endomorphisms of the
the fibre functor. Its elements are collections $\{ f_V\}$
with $f_V\in \End(\cF(V))$ such that for any morphism
$\psi: V\ra W$, one has  $\cF(\psi)f_V=f_W\cF(\psi)$.
Let $\aut^\otimes (\cF)$ be the set of elements $\{f_V\}$ of  $\End(\cF)$ 
satisfying 
$$
f_1=1 \;\;\; \;\; f_{V\otimes W}=f_V\otimes f_W.
$$
The existence of duals in $\cC$ implies that any element in 
$\aut^\otimes (\cF)$  consists entirely of automorphisms, and hence
there is no need to include a condition for invertibility.
The algebra $\End(\cF)$ is a projective limit of finite dimensional
algebras and it is endowed with a projective limit topology.
The subset $\aut^\otimes (\cF)$  has a structure of projective limit
of algebraic varieties.

Let $G$ be a group such that $\Rep(G)$, with the functor $\cF_G$ 
defined as usual, is a Tannakian category.  There is a map form $G$ to 
$\aut^\otimes (\cF_G)$ 
which sends an element $g\in G$ to the natural automorphims $\{f_V\}$
of $\cF_G$ defined by setting $f_V$ equal to the action of $g$ on
the vector space $\cF_G(V)$ underlying the representation $V$.
As mentioned above, for complex affine group schemes this map is
an isomorphism. 

We come now to the Tannakian nature of twisted Higgs bundles.

\begin{prop} \label{higgs-tannaka}
The tensor category of polystable $L$-twisted Higgs bundles $E$
 over $X$, satisfying  $\C_1(E).[\omega]^{\dim X - 1}=0$ and 
$\ch_2(E).[\omega]^{\dim X - 2}=0$,
with fibre functor  defined by sending an $L$-twisted
Higgs bundle to  the fibre of the bundle at a fixed point of $X$, is a 
neutral Tannakian category. 
\end{prop}
\pf
Let $(E,\theta)$ and $(F,\eta)$ be two $L$-twisted Higgs bundles. 
Its tensor product is given by
the $L$-twisted Higgs bundle $(E\otimes F,\theta\otimes 1 + 1\otimes
\eta)$.  The polystability of the tensor product can be proved
directly,
but it follows also from the existence of metrics satisfying 
(\ref{twisted-he}).
The tensor product of the two metrics satisfies (\ref{twisted-he}) as well and
hence the tensor product Higgs bundle is polystable by Theorem \ref{hk}.
This operation defines an associative and commutative tensor category.
Suppose $f: (E,\theta ) \ra (F, \eta )$ is a morphism of $L$-twisted
 Higgs bundles, namely a sheaf homomorphism $f:E\to F$ such that the 
appropriate diagram commutes. Then the subsheaf $V$ of $F$ 
generated by the image is
invariant under $\eta $, and the kernel of $f$ is invariant under 
$\theta$, since their generic fibres are clearly invariant.  Now we
have 
$\deg E/\ker f \geq 0$ since $E$ is semistable and $\deg V \leq 0$ 
since $F$ is semistable. But then we have a generic isomorphism
induced by $f$ from
$E/\ker f \ra V$. This is only possible if the degrees of both these
bundles are 0, and the above map is actually an isomorphism. This shows
that semistable $L$-twisted Higgs bundles of degree 0, form an abelian category.
The dual of a pair 
$(E,\theta)$ 
is the pair $(E^\ast, {\theta}^{\vee})$, where $E^\ast$ is the dual 
bundle to $E$  and 
$\theta^{\vee}$ is the map obtained by transposing $\theta $ and 
tensoring with the canonical line bundle. Obviously, the Higgs bundle 
$(\cO,0)$ is a unit and satisfies that $\End((\cO,0))=\bC$.

The fibre functor $\cF$ is defined by choosing a point $x\in X$ and 
sending $(E,\theta)$ to the fibre of $E$ at the point $x$. 
The faithfulness of the $\cF$ follows from the polystability
of $(E,\theta)$.

In this paper we shall address the problem of describing
the corresponding Tannaka group when  $\deg L=0$. We will
give an answer to this problem when the order of $L$ is finite, that is when 
$L^n$ is the trivial line bundle for some $n$.
This answer is given in terms of the Tannaka group of the category
of ordinary Higgs bundles, i.e. those for which  $L$ is the trivial line bundle.

%%%%%%%%%%%%%%%%%%%%%%%%%%%%%%%%%%%%%%%%%%%% 
\section{Ordinary Higgs bundles}\label{higgs}
%%%%%%%%%%%%%%%%%%%%%%%%%%%%%%%%%%%%%%%%%%%%
Our goal in this section is to describe the Tannaka group of the
category of polystable (untwisted) Higgs bundles $E$
a compact \kahler\ manifold $(X,\omega)$ satisfying
\be
\C_1(E).[\omega]^{\dim X - 1}=0\;\;\; \mbox{and} \;\;\;
\ch_2(E).[\omega]^{\dim X - 2}=0.
\label{vanishing}
\ee

This is done by means of the correspondence between polystable Higgs
bundles satisfying (\ref{vanishing})
and semisimple complex representations of the
fundamental group.
In the sequel we briefly recall the main ideas of this correspondence
 (see \cite{Si1,H1,D2,C} for details).

%%%%%%%%%%%%%%%%%%%%%%%%%%%%%%%%%%%%%%%%%%%%%%%%%%%%%%%%%%%%%
\subsection{Higgs bundles, flat bundles and representations of the
fundamental group}\label{higgs}
%%%%%%%%%%%%%%%%%%%%%%%%%%%%%%%%%%%%%%%%%%%%%%%%%%%%%%%%%%%%
To associate a complex representation to a Higgs bundle we will pass 
through the intermediate category of flat bundles.
Let $V$ be $C^\infty$ 
complex vector bundle  of rank $r$  over $X$, and  let $D$ be a $\GL(r,\bC)$ 
connection on $V$. We say that $D$ is {\em flat} if its curvature 
vanishes,
i.e. $D^2=0$. If $D$ is a flat connection on $V$  the pair  $(V,D)$ is
called a {\em flat bundle} since, by using the flat connection, 
one can find an open cover of $X$ with constant transition functions for $V$.
If  $D$ is a connection on a vector bundle $V$ and $x$ is a 
(fixed) point of $X$ recall that the {\em holonomy group} of $D$ is the
group of endomorphisms of $V_x$ obtained by parallel transport 
along all closed curves starting at $x$. If $D$ is flat the parallel
displacement depends only on the homotopy class of the closed curve
and defines a homomorphism
$$
\rho:\pi_1(X,x) \lra \GL(V_x),
$$
whose image is the holonomy of $D$.
Conversely, given a representation $\rho:\pi_1(X,x) \ra \GL(r,\bC)$
one can construct a vector bundle $V$ of rank $r$ with a flat connection
by setting
$$
V=\tilde{X}\times_\rho \bC^r,
$$
where $\tilde{X}$ is  the universal cover of $X$ and 
$\tilde{X}\times_\rho \bC^r$ is  the quotient of 
$\tilde{X}\times \bC^r$ by the action of $\pi_1(X,x)$ given by 
$(y,v)\mapsto (\gamma(y), \rho(\gamma)v)$ for 
$\gamma\in \pi_1(X,x)$ (regarded as the covering 
transformation group acting on $\tilde X$). The trivial connection
on $\tilde{X}\times \bC^r$ descends to give a flat connection on 
$V$, whose holonomy is the image of $\rho$.

The relation between flat bundles and Higgs bundles involves
a certain class of metrics over a flat bundle---the so-called harmonic
metrics. Let  $(V,D)$ be a flat bundle over $X$. Given a  metric $h$ on 
$V$  we can decompose $D$ uniquely as $D=\nabla+\Psi$ where $\nabla$ is
a unitary connection on $V$ and $\Psi$ is a is a 1-form with values 
in the self-adjoint endomorphisms of $V$. The metric $h$ is said 
to be {\em harmonic} if
$$
\nabla^\ast \Psi=0,
$$
where we use the  metric on $X$ to define $\nabla^\ast$.
A metric $h$ on $V$  is just a section of a certain 
$\GL(r,\bC)/U(r)$-bundle over $X$. This can be viewed as a 
$\pi_1(X)$-equivariant function 
$$
\tilde{h}:\tilde{X}\ra \GL(r,\bC)/U(r),
$$
where $\tilde X$ is the universal cover of $X$.
It turns out that $\nabla^\ast \Psi=0$
is equivalent to the condition that the map 
$\tilde h$ should be harmonic. In fact the one-form $\Psi$ can be 
identified with the differential of $\tilde h$, and $\nabla$ with
the pull-back of the Levi--Civita connection on $\GL(r,\bC)/U(r)$.

To state an existence theorem for such metrics we need the following 
definitions. 
A flat bundle $(V,D)$
is said to be {\em irreducible} if $V$ has no non-trivial $D$-invariant 
subbundles. It will be called {\em semisimple} if any $D$-invariant 
subbundle has a $D$-invariant complement. Any semisimple connection is
a direct sum of irreducible ones. 
\begin{thm}\label{don-cor} 
A flat bundle $(V,D)$ over $X$ admits a 
harmonic metric if and only if it is semisimple.
\end{thm}
This theorem is proved by Donaldson \cite{D2} for rank 2 bundles when $X$ is a Riemann surface, and
in full generality (including the base manifold being  a compact
Riemannian manifold of arbitrary dimension) by Corlette \cite{C}.

Let $(E,\varphi)$ be a Higgs bundle over $X$. 
We want to associate to it a flat bundle over $X$. 
This is not always possible, but if 
$(E,\varphi)$ supports a hermitian metric $h$ satisfying 
$\Lambda F_h+\Lambda [\varphi,\varphi^\ast]=0$ then we can consider the
pair $(V,D)$, taking $V$ to be  the underlying $C^\infty$ bundle 
to $E$ and $D=\dbar_E+\partial_h + \varphi +\varphi^\ast$, where 
$\partial_h$  is a differential operator such that 
$\dbar_E +\partial_h$ is 
the unique connection compatible with the metric and the 
holomorphic structure of $E$. A simple computation shows the 
following.
\begin{lemma}
$\dbar_E\varphi=0$, $\varphi\wedge\varphi=0$ and $F_h +[\varphi,\varphi^\ast]=0$ imply that
 $D$ is flat, i.e. $D^2=0$.
\end{lemma}

We will show now how to associate a Higgs bundle to a flat bundle.
Let $(V,D)$ be a flat bundle over $X$.
We want to produce out of it a stable Higgs bundle $(E,\varphi)$  
over $X$. Let $h$ be a Hermitian
metric on $V$. We can  decompose $D$ in its $(1,0)$ and $(0,1)$ 
componentes
$$
D=D'+D''
$$
and consider the unique operators $D_h''$ and $D_h'$ so that
$D'+D_h''$ and $D_h'+D''$ become $h$-unitary connections.
Let
$$
\partial_h=\frac{D'+D'_h}{2}\;\;\;
\dbar_h=\frac{D''+D''_h}{2}\;\;\;
\varphi_h=\frac{D'-D'_h}{2}\;\;\;
\varphi_h^\ast=\frac{D''-D''_h}{2}.
$$
It is not difficult to see that
\begin{lemma}
$D^2=0$ implies that $\varphi_h\wedge\varphi_h=0$ and $F_h+[\varphi_h,\varphi_h^\ast]=0$. 
Where $F_h$ is the
curvature of $\dbar_h+\partial_h$.
\end{lemma}

Of course $\dbar_h$ defines a holomorphic structure on $V$, but
$\varphi_h$ need not be holomorphic with respect to it, i.e. 
there is 
no reason why $\dbar_h\varphi_h=0$. 
This happens precisely when the metric is harmonic (\cite{D2,H1,Si2}).

Putting everything together one has the following.
\begin{thm}\label{3-tannaka}
There is a tensor functor which is an equivalence of categories
between the category of polystable Higgs bundles over $(X,\omega)$ satisfying
(\ref{vanishing}) and the
category of semisimple flat bundles which, in turn, is 
equivalent to the tensor category of semisimple complex
representations of the fundamental group of $X$.
\end{thm}

%%%%%%%%%%%%%%%%%%%%%%%%%%%%%%%%%%%%%%%%%%%%%%%%%%%%%%%%
\subsection{The Tannaka group for ordinary Higgs bundles}
%%%%%%%%%%%%%%%%%%%%%%%%%%%%%%%%%%%%%%%%%%%%%%%%%%%%%%%%

One can easily prove the following \cite[Lemma 6.1]{Si2}.

\begin{prop}\label{discrete}
Let $H$ be a finitely generated group. The tensor category of
semisimple representations of $H$, with its obvious fibre functor, is a neutral
Tannakian category whose Tannaka group is naturally isomorphic to
the pro-reductive completion of $H$. 
\end{prop}
We recall that the pro-reductive completion 
of a group $H$ is 
a projective limit $ G={\underline \lim} (\Gamma,\rho)$,
where the inverse limit runs over the directed system of representations 
$\rho: H \ra\Gamma$ for complex reductive groups 
$\Gamma$ (we shall assume that
the image of $\rho$ is Zariski dense for convenience). An arrow 
$(\Gamma,\rho)\ra (\Gamma',\rho')$ consists of a homomorphism 
$f:\Gamma \ra \Gamma'$ such that $f\rho=\rho'$. 
The group $G$ is characterised by the following universal property: 
For every  representation $H\ra\Gamma$ into a complex reductive
group there exists a unique extension 
$G\ra\Gamma$ such that the following  diagram
$$
\begin{array}{ccc}
H & \ra & \Gamma\\
\downarrow & \nearrow &        \\
G          &          &  
\end{array}
$$
commutes.

From this proposition and Theorem 
\ref{3-tannaka} one concludes the following.

\begin{thm}
Let $G$ be the Tannaka group of the 
category of polystable Higgs bundles over $X$ satisfying
(\ref{vanishing}). 
Then $G$ is naturally isomorphic to the pro-reductive completion of $\pi_1(X,x)$.
\end{thm}

One of the main ingredients for the description of the Tannaka group
in the twisted situation is a natural action of the group $\bC^\ast$ 
on the category of polystable Higgs bundles given by  
$$
(E,\varphi)\mapsto (E,\lambda\varphi)
\;\;\;\mbox{for every}\;\; \lambda \in\bC^\ast.
$$
This action induces an action of $\bC^\ast$ on the category of semisimple
representations of the fundamental group. It should be pointed out that,
while this action is very clear and explicit from the point of view
of Higgs bundles, its explicit effect on a representation of the 
fundamental group is not easy to describe.
  
One can formalise the action of $\bC^\ast$ on a Tannakian category 
$(\cC,\cF)$ in  terms of certain tensor functors satisfying canonical 
axioms. If the action preserves the fibre functor $\cF$ one has
an action of $\bC^\ast$ on $\End(\cF)$ by sending the element 
$\{f_V\}$ of $\End(\cF)$ to $\{f^\lambda_V\}$ with 
$f^\lambda_V=f_{\lambda V}$ for every $\lambda\in \bC^\ast$,
and hence one has an action on the Tannaka group $\aut^\otimes(\cF)$.
The action of $\bC^\ast$ on the category of polystable pairs preserves
clearly the fibre functor since the bundle is unchanged, and one 
can then transfer this action to $G$---the pro-reductive completion
of the fundamental group. More precisely one has  the following theorem 
(\cite{Si2}[Theorem 6]).

\begin{thm}\label{action}
There exists a unique action of $\bC^\ast$ on $G$, each 
$\lambda\in\bC^\ast$ acting
by a homomorphism of pro-reductive groups, such that if 
$\rho:G\ra\GL(n,\bC)$ is the 
representation corresponding to $(E,\varphi)$, then 
$\rho\circ\lambda$ is the 
representation corresponding to $(E,\lambda\varphi)$.
\end{thm}

%%%%%%%%%%%%%%%%%%%%%%%
\section{Main theorem}
%%%%%%%%%%%%%%%%%%%%%%%

Before stating our main result we will prove some preliminary 
necessary facts on groups with a $\bC^\ast$-action.

%%%%%%%%%%%%%%%%%%%%%%%%%%%%%%%%%%%%%%%%%%%%%%%%%
\subsection{Twisted groups}\label{twisted-groups}
%%%%%%%%%%%%%%%%%%%%%%%%%%%%%%%%%%%%%%%%%%%%%%%%%

Let $G$ be a group  and let $\bC^\ast$ act on $G$, i.e.
we have a homomorphism
$$
\bC^\ast\lra\aut(G).
$$
For every $\lambda\in \bC^\ast$ and $g\in G$ we shall denote
by $g^\lambda$ the image of $g$ by the automorphism $G\ra G$
defined by $\lambda$.

\begin{prop}
Let $\chi: G\ra \bC^\ast$ be a character satisfying
\be
\chi(g)=\chi(g^\lambda) \;\;\; \mbox{for every}\;\; g\in G,\lambda\in\bC^\ast.
\label{star}
\ee

We can define a group $G_\chi$ by taking the underlying set to be $G$
and the group operation to be 
\be
g\twist h=gh^{\chi^{-1}(g)}\;\;\;\mbox{for every}\;\; g,h\in G.
\label{twist-operation}
\ee
\end{prop}
\pf
Associativity results  from the following computation for  $g,h,k\in G$.
\begin{eqnarray}\nonumber
g\twist(h\twist k)&=& g(h\twist k)^{\chi^{-1}(g)}\\ \nonumber
                  &=& g(hk^{\chi^{-1}(h)})^{\chi^{-1}(g)}\\ \nonumber
                  &=& gh^{\chi^{-1}(g)}k^{\chi^{-1}(h)\chi^{-1}(g)}\\\nonumber
                  &=& gh^{\chi^{-1}(g)}k^{\chi^{-1}(gh)}. \nonumber
\end{eqnarray}

\begin{eqnarray}\nonumber
(g\twist h)\twist k&=& gh^{\chi^{-1}(g)}k^{\chi^{-1}(gh^{\chi^{-1}(g)})}\\ \nonumber
                    &=& gh^{\chi^{-1}(g)}k^{\chi^{-1}(g)\chi^{-1}(h^{\chi^{-1}(g)})}\\
                    \nonumber
                    &=& gh^{\chi^{-1}(g)}k^{\chi^{-1}(g)\chi^{-1}(h)}\;\; 
                                                \mbox{by condition\ }
                    (\ref{star})\\\nonumber
                    &=& gh^{\chi^{-1}(g)}k^{\chi^{-1}(gh)}. \nonumber
 \end{eqnarray}
The identity element $e\in G$ is also the identity
element of $G_\chi$. Indeed,
$$
g\twist e= g e^{\chi^{-1}(g)}=ge=g.
$$
Let $g\in G$. Take $h$ to be the preimage of $g^{-1}$
under the automorphism $\chi^{-1}(g):G\ra G$, i.e. $h^{\chi^{-1}(g)}=g^{-1}$. Then
$$
g\twist h=gh^{\chi^{-1}(g)}=gg^{-1}=e
$$
and hence $h$ is the inverse of $g$.
We have then proved that $G_\chi=(G,\twist)$ is a group.
\qed

We shall now prove few properties about $G_\chi$ that will be useful
later.
\begin{prop}
Let $\chi:G\ra\bC^\ast$ be  a character. Then the map
$\alpha_\chi:G_\chi\ra \bC^\ast$  defined by $\alpha_\chi(g)=\chi(g)$
is a character of $G_\chi$. In particular $K=\ker \chi$ is a
normal subgroup of $G_\chi$.
\end{prop}
\pf
\begin{eqnarray}\nonumber
\alpha_\chi(g\twist h)&=& \chi(gh^{\chi^{-1}(g)})\\\nonumber
                      &=& \chi(g) \chi(h^{\chi^{-1}(g)})\\\nonumber
                      &=& \chi(g) \chi(h) \;\;\mbox{by (\ref{star})}\\\nonumber
                      &=& \alpha_\chi(g)\alpha_\chi(h) \;\; \mbox{for every}
\;\; g,h\in G_\chi.\nonumber
\end{eqnarray}
\qed
\begin{prop}
The action of $\bC^\ast$ on $G$ defines an action of $\bC^\ast$ on $G_\chi$.
Moreover, this action leaves $K$ invariant.
\end{prop}
\pf
\begin{eqnarray}\nonumber
(g\twist h)^\lambda & = & (gh^{\chi^{-1}(g)})^\lambda\\\nonumber
                    & = & g^\lambda h^{\chi^{-1}(g)\lambda}\\\nonumber
                    & = & g^\lambda h^{\lambda\chi^{-1}(g)}\\\nonumber
                    & = & g^\lambda h^{\lambda\chi^{-1}(g^\lambda)} \;\; 
                                            \mbox{by (\ref{star})}\\\nonumber
                    & = & g^\lambda \twist  h^\lambda.
\end{eqnarray}
To see that the action of $\bC^\ast$ on $G_\chi$ leaves $K$ invariant,
let $g\in K$, and $\lambda\in \bC^\ast$. By (\ref{star}), we see that
$$
\alpha_\chi(g^\lambda)=\chi(g^\lambda)=\chi(g)=e,
$$
and hence $g^\lambda \in K$.
\qed

%%%%%%%%%%%%%%%%%%%%%%%%%%%%%
\subsection{Main theorem}
%%%%%%%%%%%%%%%%%%%%%%%%%%%%%

Let $L$  be a holomorphic line bundle over $X$ such that $\deg L=0$.
Let  $\pi_1(X)=\pi_1(X,x)$ be  the 
fundamental group of $X$ with respect to a fixed point $x\in X$. 
The line bundle $L$ corresponds to a unitary character
 $\chi': \pi_1(X) \ra \U(1)$. That is, 
if $\tilde{X}$ is the universal cover of $X$, $L$  is
the line bundle associated to the $\pi_1(X)$-principal bundle 
$\tilde{X}\ra X$ via the representation $\chi'$.

Let $G$ be the pro-reductive completion of $\pi_1(X)$. As 
discussed in section \ref{tannaka}, $G$ is isomorphic to the Tannaka
group of the  category of polystable Higgs bundles and there is
 $\bC^\ast$ on $G$, induced from 
the action of $\lambda\in \bC^\ast$ on a Higgs bundle $(E,\varphi)$ given by
$(E,\varphi)\mapsto (E,\lambda\varphi)$ (Theorem \ref{action}).
We shall apply the construction of a twisted group structure given above
to the pro-reductive completion of the fundamental
group of $X$.

\begin{prop}\label{twisted}
Let $\chi:G\ra\bC^\ast$ be the extension to $G$ of a {\em unitary}
character $\chi':\pi_1(X)\ra \U(1)$. Then the action of $\bC^\ast$ on $G$ 
considered above satisfies (\ref{star}), i.e.
$$
\chi(g)=\chi(g^\lambda)\;\;\;
\mbox{for every}\;\; g\in G, \lambda \in\bC^\ast.
$$
\end{prop}
\pf
We shall regard $G$ as the group of tensor automorphisms of the fibre 
functor of the Tannakian category of Higgs bundles over $X$ 
(see section \ref{tannaka}). Hence
$g\in G$ basically associates  to any stable Higgs  pair $(E,\varphi)$ 
an automorphism $f_{(E,\varphi)}$ of the fibre $E_x$, in a functorial
way. Let  $\lambda\in \bC^\ast$. The element 
$g^\lambda\in G$ associates to $(E,\varphi)$ the automorphism 
$f_{(E,\lambda\varphi)}$.

Let now $\rho$ be a semisimple representation  of $G$. Of course $\rho$
gives rise to a semisimple representation of $\pi_1(X,x)$ and hence to a 
polystable Higgs bundle $(E,\varphi)$. One can see that 
$$
\rho(g)=f_{(E,\varphi)}.
$$
But if $\rho|_{\pi_1(X)}$ is unitary, then the associated Higgs bundle is
$(E,0)$, and hence $\rho(g)=\rho(g^\lambda)$. In particular, for our 
``unitary'' character we have $\chi(g)=\chi(g^\lambda)$, which concludes 
the proof.
\qed

We can thus consider the group $G_\chi$ 
by means of the construction given in the previous  section.

We are now ready to state the main theorem of this paper.
\begin{thm}\label{main-theorem}
Let $L$ be a holomorphic line bundle of degree zero and finite order.
Then $G_\chi$ is the Tannaka group of the category of 
polystable $L$-twisted Higgs bundles satisfying (\ref{vanishing}).
More precisely, there is a tensor functor which is an equivalence of
categories between the
category of  polystable $L$-twisted Higgs bundles 
satisfying (\ref{vanishing})
and the category of semisimple complex representations of $G_\chi$.
\end{thm}

It is natural to conjecture that this result is also true in the 
infinite order case. We hope to come back to this in a 
future paper.

%%%%%%%%%%%%%%%%%%%%%%%%%%%%%%%%%%%%%%%%%%%%%%%%%%%%%%%%%%%%%%%%%%%%%%
\section{Proof of main theorem}
\label{finite}
%%%%%%%%%%%%%%%%%%%%%%%%%%%%%%%%%%%%%%%%%%%%%%%%%%%%%%%%%%%%%%%%%%%%%%
Let $\chi':\pi_1(X)\ra \U(1)$ be the unitary character corresponding 
to $L$ and let $\Gamma=\im\chi'$ be the image.
Let $\tilde{X}$ be the universal cover of $X$ and let 
$Y=\tilde{X}/\Ker \chi'$. Then $\Ker \chi'\cong \pi_1(Y)$ and 
one has the covering map
$$
p:  Y\lra X
$$
whose Galois group is $\Gamma=\im\chi'\cong \pi_1(X)/\pi_1(Y)$.

The basic strategy that we shall follow in our approach is to translate
our  problem into a problem on  $Y$:  
If the line bundle $L$ is of finite order, i.e. if  a finite power of
 $L$ is isomorphic to the trivial line bundle, $\Gamma$ is a finite
cyclic group and hence $Y$ is compact. 
The pull-back of $L$ to $Y$ is of course trivial, and we can thus 
use the untwisted theory over $Y$ to prove our theorem. 

%%%%%%%%%%%%%%%%%%%%%%%%%%%%%%%%%%%%%%%%%%%%%%%%%%%%%%%%%%%%%%%%%%%%
\subsection{Twisted Higgs bundles and representations of $\pi_1(Y)$}
\label{twisted-N}
%%%%%%%%%%%%%%%%%%%%%%%%%%%%%%%%%%%%%%%%%%%%%%%%%%%%%%%%%%%%%%%%%%%%
Let $(E,\theta)$ be an $L$-twisted Higgs bundle over $X$. 
Its pull-back to $Y$, $(F,\varphi)=(p^\ast E, p^\ast \theta)$, 
becomes a genuine Higgs bundle. It is not difficult to 
characterize  Higgs bundles 
on $Y$ that {\em come} from twisted 
Higgs bundles on $X$.

\begin{prop} 
There is an equivalence of tensor categories  between 
the category of $L$-twisted Higgs bundles $(E,\theta)$ on $X$ and 
the category of Higgs 
bundles $(F,\varphi)$ on $Y$ that satisfy 
\be
\gamma^\ast F=F\;\;\;\mbox{and}
\;\;\gamma^\ast\varphi=\gamma\varphi\;\;
\mbox{for every}\;\; \gamma\in \Gamma.
\label{h-act}
\ee
\end{prop}
\pf
We must first clarify our notation: On the one hand
we  regard
$\Gamma$ as the Galois group of the cover
$Y\ra X$, and on the other as a subgroup of $\U(1)$.
So when we write $\gamma^\ast\varphi$ we are thinking of $\gamma$
as a transformation of $Y$, while when we write $\gamma\varphi$
we regard $\gamma$ as an element of $\U(1)$ that multiplies the 
Higgs field.

One direction is clear:  If $(E,\theta)$ is  an $L$-twisted 
Higgs bundle over $X$ its pull-back, 
$(F,\varphi)=(p^\ast E,p^\ast\theta)$, obviously satisfies (\ref{h-act}).
To see the converse, observe that the first condition in (\ref{h-act}) 
amounts to saying 
that the bundle $F\ra Y$ descends to a bundle $E\ra X$, since the Galois group is cyclic. 
It is obvious that this correspondence is compatible with tensor products.

\begin{prop}
Let $(E,\theta)$ be an $L$-twisted Higgs bundle on $X$, and let 
$(F,\varphi)$ its pull-back to $Y$.
The pair $(E,\theta)$ is polystable if and only if
$(F,\varphi)$ is polystable.
\end{prop}
\pf
The vector bundle $F=p^\ast E$ is a $\Gamma$-equivariant vector bundle,
and one can consider for the pair $(F,\varphi)$ a weaker stability
notion  consisting of  
the usual  stability condition, but
only for $\Gamma$-equivariant subsheaves of $F$. One can show \cite{G}
that this $\Gamma$-equivariant 
condition for polystability is equivalent to the 
polystability of $(F,\varphi)$. Now, suppose that $(E,\theta)$ is not
stable. The pull-back of the destabilizing subsheaf would violate
the $\Gamma$-equivariant stability of $(F,\varphi)$.
Conversely, if $(F,\varphi)$ is not $\Gamma$-equivariantly stable, 
then let $F'\subset F$ be a $\Gamma$-equivariant $\varphi$-invariant
destabilizing subbundle. Because of $\Gamma$-equivariance $F'$
descends to a $\theta$-invariant destabilizing subbundle of $E$.

Another way of proving this proposition is to use 
Theorem  \ref{hk}: If $(F,\varphi)$ is polystable there exists 
a Hermitian metric on $F$ solving equation (\ref{twisted-he}). By the uniqueness of 
the solution, this metric must be $\Gamma$-invariant and hence 
descend to a metric on $E$ solving (\ref{twisted-he}) implying the 
polystability of $(E,\theta)$. The converse is also clear.
\qed

Since, as we know, there is a tensor functor which is an equivalence of
categories between the category of polystable Higgs bundles over $Y$ 
and the category of semisimple representations of the fundamental group 
of $Y$, 
we shall now give an interpretation of condition (\ref{h-act}) 
in terms of representations of $\pi_1(Y)$.
To do this we need to understand the two actions of $\Gamma$ 
on Higgs bundles in terms of the representations of $\pi_1(Y)$. 
To explain the action corresponding
to pull-backing  by $\gamma\in\Gamma$, we shall digress 
 a little.

Let $G$ be a group and let $K\subset G$ be a normal subgroup.
Let $\Gamma=G/K$. We have a short exact sequence
$$
1\lra K \lra G \lra \Gamma\lra 1.
$$
The group $G$ acts on the set of representations of $K$, via inner 
automorphisms, i.e. if $\rho$ is a representaion of $K$,  $g\in G$
sends $\rho$ to $\rho\circ \Int_g|_K$, where 
$$
\Int_g|_K (h)= g h g^{-1}\;\;\mbox{for every}\;\; h\in K.
$$
\begin{prop}\label{gamma-action}
Let $\gamma\in \Gamma$ and let $g_\gamma\in G$ be a lift of $\gamma$.
The map 
\be
 [\rho]\mapsto  \gamma\cdot [\rho]= [\rho\circ\Int_{g_\gamma}|_K]
\;\;\;\mbox{for every}\;\;\;[\rho]\in\Rep(K)\;\;\mbox{and}\;\;
\gamma\in \Gamma \label{action-rep}
\ee
defines
an action of $\Gamma$ on the set $\Rep(K)$  of equivalence classes of
representations of  $K$.
\end{prop}
\pf 
It is clear since two lifts of $\gamma$ differ by an element of $K$.
\qed

The following is immediate.
\begin{prop}\label{equivalence-action}
 Let $\rho$ be a representation of $K$ and let 
$\gamma\in \Gamma$. Then
$\gamma\cdot [\rho]=[\rho]\; \mbox{for every}\; \gamma\in \Gamma$
is equivalent to
$[\rho]=[\rho\circ \Int_g]\;\;\mbox{for every}\;g\in G$.
\end{prop}

In our situation we  have the extension
$$
1\lra \pi_1(Y)\lra \pi_1(X) \lra \Gamma\lra 1.
$$
By the Proposition   \ref{gamma-action}  there is
an action of $\Gamma$ on 
representations of  $\pi_1(Y)$, given for every 
$[\rho]\in \Rep(\pi_1(Y))$ and 
$\gamma\in \Gamma$ by
$[\rho]\mapsto\gamma\cdot [\rho]= [\rho\circ\Int_{g_\gamma}]$,
where $g_\gamma\in $ is any lift of $\gamma$ to $\pi_1(X)$.
The following is clear.
\begin{prop}
The action of $\Gamma$ on $\Rep(\pi_1(Y))$ given above corresponds 
to the action of $\Gamma$ on the set of equivalence classes of Higgs 
bundles  and flat bundles over $Y$ defined by 
$\gamma\cdot (F,\varphi)=(\gamma^\ast F,\gamma^\ast \varphi)$
and $\gamma\cdot (V,D)=(\gamma^\ast V,\gamma^\ast D)$, respectively.
\end{prop}

As we saw in section \ref{tannaka}, the action of $\bC^\ast$ on the  
moduli of stable Higgs bundles on
$Y$ given for $\lambda \in\bC^\ast$, by 
$[(F,\varphi)]^\lambda=[(F,\lambda \varphi)]$ defines an 
action on $\Rep(\pi_1(Y))$. If $[\rho]$ corresponds to $[(F,\varphi)]$ we
denote by $[\rho]^\lambda$ the representation corresponding
to $[(F,\lambda \varphi)]$.
The following proposition follows immediately.

\begin{prop}
Let $(F,\varphi)$ be a polystable Higgs bundle over $Y$,
and let $[\rho]$ be the corresponding semisimple representation of 
$\pi_1(Y)$. 
The  condition (\ref{h-act}) is equivalent to
\be
\gamma\cdot [\rho]=[\rho]^\gamma \;\;\mbox {for every}\;\;\gamma\in
\Gamma.
\label{r-act}
\ee
\end{prop}
We would like now to relate the representations of $\pi_1(Y)$ satisfying
(\ref{r-act}) to representations of $G_\chi$. To do this we will relate 
them first to the representations of $K=\Ker \chi$. 

%%%%%%%%%%%%%%%%%%%%%%%%%%%%%%%%%%%%%%%%%%%%%%%%%%%%%%%%%%%%%%%%%%%%%%%%
\subsection{From representations of $\pi_1(Y)$ to representations of $K$}
%%%%%%%%%%%%%%%%%%%%%%%%%%%%%%%%%%%%%%%%%%%%%%%%%%%%%%%%%%%%%%%%%%%%%%%%%
\begin{prop} 
If $\chi':\pi_1(X)\ra \U(1)$ be a unitary character of finite order of
$\pi_1(X)$. Let $\chi$ be its extension to $G$, the Tannaka closure
of $\pi_1(X)$, and let $K=\Ker \chi$. Then any
complex representation of $\pi_1(Y)$ extends to a complex representation
of $K$.
\end{prop}
\pf
Let $n$ be the order of $\chi'$. Since ${\chi'}^n$ is trivial, we get from the uniqueness
of extensions that $\chi^n$ is also trivial.
Hence we have the diagram
$$
\begin{array}{ccccccccc}
1  & \lra &  K    & \lra &  G &  \lra &  \Gamma & \lra &  1\\
   &      & \uparrow &   & \uparrow &  & ||     &      &   \\
1  & \lra & \pi_1(Y) & \lra & \pi_1(X) &  \lra  & \Gamma & \lra & 1.
\end{array}
$$
Let $(V,\rho)$ be an irreducible representation of $\pi_1(Y)$. 
This induces
a semisimple representation $(W,\ind (\rho))$ of $\pi_1(X)$, where
$$
W=\{f:\pi_1(X)\ra V\;\;|\;\; f(hy)=\rho(y^{-1}) f(h),\;\;\mbox{for every}\; 
y\in \pi_1(Y),\; h\in \pi_1(X)\}.
$$
The morphism $\ind(\rho)\in \aut (W)$ is defined by
$$
(hf)(x)=\ind (\rho)_h(f)(x)= f(h^{-1}x)\;\;\;
\mbox{for every}\; h,x\in \pi_1(X)\; \mbox{and} \; f\in W.
$$
This representation extends of course to a representation of $G$ since
$G$ is the Tannaka closure of $\pi_1(X)$. Let $\cF(G)$ be the set of 
polynomial functions on $G$ with values in $V$, and let $e:W\ra V$ be the 
evaluation map $f\mapsto f(1)$.
Consider the map $\Phi:W\ra \cF(G)$, defined by
$$
\Phi(w)(g)=e(g^{-1}w).
$$
Let us also consider the restriction map 
$\Res:\cF(G)\ra \cF(\pi_1(X))$, where 
$\cF(\pi_1(X))$ is the set of functions on $\pi_1(X)$. The image of 
$\Res$ is contained in $W$ and  the composite 
$W\stackrel{\Phi}{\lra}\cF(G)\stackrel{\Res}{\lra} W$ is the identity
since $\Phi(w)(h)=e(h^{-1}w)= (h^{-1}w)(1)=w(h)$, that is 
$\Phi(w)|_{\pi_1(X)}=w$.

The affine group $G$ is defined as Spec of the algebra of representation
functions (namely coefficients of semisimple representations) of $\pi_1(X)$.
So the restriction map of the algebra of functions on $G$ to the Zariski
closure of $\pi_1(X)$ is onto (because it is affine) and one-to-one,
and hence they are the same. In other words $\pi_1(X)$ is dense in $G$.
Now, the image of $W$ in $\cF(G)$
satisfies
$$
f(gy)=\rho(y^{-1}) f(h),\;\;\mbox{for every}\; 
y\in \pi_1(Y),\; g\in G.
$$
This is because $\psi(g)=f(gy)-\rho(y^{-1})f(g)$ is 0
on $\pi_1(X)$ and the image of  $\pi_1(X)$ is Zariski dense on $G$. 
Let 
$$
W'=\{f\in \cF(G)\;\;|\;\; f(gy)=\rho(y^{-1}) f(g),\;\;\mbox{for every}\; 
y\in \overline{\pi_1(Y)},\; g\in G\}.
$$
The restriction map $W'\ra \cF(\pi_1(X))$ is injective and the image is 
$W$,
therefore $W'=W$ and hence $G$ acts on it and the dimension of $W$ must
be $|G/\overline{\pi_1(Y)}|\dim V$ which has to coincide 
with $|G/K|\dim V$.
Hence $\overline{\pi_1(Y)}=K$.

We shall now show  that $K$ acts on $\Ker e$ and hence on $V$. 
Let $f,f'\in W$. From $(f-f')(1)=0$ we have $(f-f')(y)=0$ for every
$y\in \pi_1(Y)$ and since $\overline{\pi_1(Y)}=K$ we have 
$(f-f')(k)= \rho(k) (f-f')(1)=0$, which complets the proof.
\qed

\begin{cor} 
$K$ is the  Tannaka closure of $\pi_1(Y)$.
\end{cor}

We thus conclude 
that if $[\rho]$ is a representation of $\pi_1(Y)$
satisfying (\ref{r-act}) it extends to a representation of $K$
satisfying the same condition, where the first action of $\Gamma$ 
on this representation is defined via the extension
$$
1\lra K \lra G \lra \Gamma\lra 1,
$$
and (\ref{action-rep}).
On the other hand, since $K$ is the Tannaka closure of $\pi_1(Y)$, 
by Theorem \ref{action}, there is an action of $\bC^\ast$ on $K$  such
that
$$
[\rho]^\lambda=[\rho\circ\lambda]
\;\;\;\mbox{for every}\;\; \lambda\in\bC^\ast.
$$
Recall that we are identifying $\lambda\in\bC^\ast$ with the
homomorphism $K\ra K$ that it defines.

Combining all this with the results of the previous subsection, we have 
the  following.
\begin{prop}
There is a functor which is an equivalence of categories between 
the category of polystable $L$-twisted Higgs bundles over $X$ and 
the category of semisimple representations $[\rho]$ of $K$ satisfying
\be
\gamma\cdot [\rho]=[\rho]^\gamma \;\;\mbox {for every}\;\;\gamma\in
\Gamma.
\label{r-act-2}
\ee
\end{prop}
From section \ref{twisted-groups}, we know that $K$ is also a 
subgroup of $G_\chi$. We shall show now that
the representations of $K$ corresponding to twisted Higgs bundles extend
to representations of $G_\chi$. To see this we shall briefly analyse
separately the general problem of extending a 
representation of a normal subgroup to the whole group.

%%%%%%%%%%%%%%%%%%%%%%%%%%%%%%%%%%%%%%%%%%%%%%%%%%%%%%%%%%%
\subsection{Extending representations of a normal subgroup}
%%%%%%%%%%%%%%%%%%%%%%%%%%%%%%%%%%%%%%%%%%% %%%%%%%%%%%%%%%%

\begin{prop} \label{ext-0}
Let $G$ be a group and  $K\subset G$ be a normal subgroup.
Let $(V,\rho)$ be a representation of $K$
which extends to a representation of $G$, then for every 
$g\in G$, the representations  $\rho$  and $\rho\circ \Int_g$ are equivalent.
\end{prop}
\pf 
Suppose that $(V,\rho)$ extends to a representation
$(V,\tilde\rho)$ of $G$. Then for every $g\in G$ and  $h\in K$ we have
$$
(\rho\circ \Int_g)(h)=\rho(ghg^{-1})=\tilde\rho(ghg^{-1})=
\tilde\rho(g) \tilde\rho(h) \tilde\rho(g^{-1}) =
\tilde\rho(g) \rho(h) \tilde\rho(g^{-1}).
$$
Thus $\rho$ and $\rho\circ \Int_g$ are equivalent.
\qed

We are interested in the case in which $\Gamma=G/K$ is a finite cyclic
group. In this situation, the converse of the previous Proposition
is actually true. More precisely:

\begin{prop}\label{cyclic}
Let $K$ be a normal subgroup of $G$ and let $\Gamma=G/K$ be a finite 
cyclic group. Let $(V,\rho)$ be a semisimple representation of $K$ such that 
$\gamma\cdot [\rho]=[\rho]$ for every $\gamma\in \Gamma$. Then $(V,\rho)$
extends to a semisimple representation of $G$.
\end{prop}
\pf
Assume first that $\rho$ is irreducible.
Let $n$ be the order of $\Gamma$ and let $a$ be a generator of $\Gamma$.
Let  $g\in G$ be a lift of $a$. 
By assumption 
there exists a matrix $T\in\GL(V)$ such that 
\be
 \rho(gkg^{-1})=T\rho(k) T^{-1} \;\;\;\mbox{for every}\; k\in K. 
\label{equivalence}
\ee
This implies that $T^n\rho(k)T^{-n}=\rho(g^nkg^{-n})$. But $g^n\in K$
and so we have  $T^n\rho(k)T^{-n}=\rho(g^n)\rho(k)\rho(g^{-n})$. In
other words, $\rho(g^{-n})T^n$ commutes with $\rho(k)$ for all $k\in
K$. Since $\rho$ is irreducible,  $\rho(g^{-n})T^n=\lambda I$
for some constant $\lambda\in\bC$, 
by Schur's Lemma, 
and replacing $T$ by $\mu T$ with $\mu^n=\lambda$, we see that 
$T^n=\rho(g^n)$. We will now extend $\rho$ to a representation
$\tilde\rho$ of $G$ by setting $\tilde\rho (g)=T$. Since $G$ is
defined as the quotient of the free product of $K$ and $\Gamma$ with
the relation $g^n=k$, for some $k\in K$, we see that $\tilde\rho$ is indeed a
representation of $G$ extending $\rho$.

Suppose now that $\rho=\bigoplus_{i=1}^m \rho_i$, where the $\rho_i's$
are irreducible, inequivalent representations such that 
\be
\rho_i=\rho_1\circ\Int_{g^{i-1}}|_K\;\;\; \mbox{for}\;\;\; 1\leq i\leq m,
\label{cyclic-rep}
\ee
with $m$ dividing $n$. Every 
representation satisfying the hypothesis of the theorem is a direct
sum of representations of this kind, and it will hence be enough to
consider this case.

The automorphism $T$ in (\ref{equivalence}) has a block decomposition
of the form

$$
T=\left(\begin{array}{cccc}
0&...& 0& A_m\\
A_1 &0& ...& 0\\
...&...& ... & ...\\
0&...&A_{m-1} & 0
\end{array}
\right),
$$
where $A_i$ with $1\leq i\leq m$ is an invertible linear transformation.
In terms of this decomposition (\ref{equivalence}) is equivalent to
\be
\begin{array}{ccc}
A_m\rho_m(k)A_m^{-1}& = & \rho_1(gkg^{-1})\\ 
A_1\rho_1(k)A_1^{-1}& = & \rho_2(gkg^{-1})\\ 
...                 & ... & ...           \\ 
A_{m-1} \rho_{m-1}(k) A_{m-1}^{-1}&=& \rho_m(gkg^{-1}). 
\end{array} \label{equivalence-2}
\ee
This implies that
$$
\rho_1(k)=A_m
A_{m-1}...A_1\rho_1(g^mkg^{-m})A_1^{-1}...A_{m-1}^{-1}A_m^{-1},
$$ 
and similarly for the other $\rho_i's$. More precisely, let
$$
B_i=\prod_{j=1}^{i-1}A_{i-j}\prod_{l=0}^{m-i}A_{m-l}.
$$
We have 
$$
\rho_i(k)=B_i\rho_i(g^mkg^{-m})B_i^{-1}
\;\;\mbox{for}\;\; 1\leq i\leq m.
$$
Iterating this, we obtain
\be
\rho_i(k)=B_i^p\rho_i(g^nkg^{-n})B_i^{-p}
\;\;\mbox{for}\;\; 1\leq i\leq m, \label{equivalence-3},
\ee
where $p=n/m$. Since $g^n\in K$,
$$
\rho_i(k)=B_i^p\rho_i(g^n)\rho_i(k)\rho_i(g^{-n})B_i^{-p}
\;\;\mbox{for}\;\; 1\leq i\leq m,
$$
from which, by Schur's Lemma, we obtain
\be
\rho_i(g^n)=\lambda_i B_i^p
\;\;\mbox{for}\;\; 1\leq i\leq m \label{equivalence-4}
\ee
for constants $\lambda_i\in\bC$.
Now, taking $k=g^n$ in (\ref{equivalence-2}), we get
\be
\tr \rho_i(g^n)=\tr \rho_j(g^n) \;\;\mbox{for}\;\; 1\leq i,j\leq m.
\label{trace}
\ee
But $\tr (B_i^p)=\tr (B_j^p)$, for $ 1\leq i,j\leq m$, since 
$A_iB_i^pA_i^{-1}=B_{i+1}^p$, and hence (\ref{equivalence-4}) and (\ref{trace})
imply that $\lambda_i=\lambda_j=\lambda$ for $ 1\leq i,j\leq m$.
As above, we can replace $T$ by $\mu T$ with $\mu^n=\lambda$, 
so that $T^n=\rho(g^n)$, and  extend $\rho$ to a representation
$\tilde\rho$ of $G$ by setting $\tilde\rho (g)=T$. 

The proof is now complete since, as mentioned above, every
representation satisfying the hypothesis of our theorem is a direct sum of
representations of the kind defined by (\ref{cyclic-rep}), and T can be
decomposed in diagonal blocks to which we can apply the above 
argument.
\qed

%%%%%%%%%%%%%%%%%%%%%%%%%%%%%%%%%%%%%%%%%%%%%%%%%%%%%%%%%%%%%%%%%%%%%%%
\subsection{From representations of $K$ to representations of $G_\chi$}
%%%%%%%%%%%%%%%%%%%%%%%%%%%%%%%%%%%%%%%%%%%%%%%%%%%%%%%%%%%%%%%%%%%%%%%
Coming back to our main theme,
recall that $K=\ker \chi$ is a normal subgroup of $G_\chi$, and we
want to see that
the representations of $K$ coming from twisted Higgs bundles extend
to representations of $G_\chi$. More precisely:
\begin{prop}
There is a tensor functor which is an equivalence of categories between 
the category of semisimple
representations of $K$ satisfying (\ref{r-act-2})
and the category of semisimple representations of $G_\chi$.
\end{prop}
\pf
Let us  consider the extension
$$
1\lra K \lra G_\chi \lra \Gamma_\chi\lra 1,
$$
where $\Gamma_\chi=G_\chi/K$.
By Proposition \ref{gamma-action} $\Gamma_\chi$ acts on $\Rep(K)$. 
Let us denote 
this action by
$\gamma\cdot_\chi [\rho]$, for every $[\rho]\in\Rep(K)$ and 
$\gamma\in\Gamma_\chi$. 
We will show now that (\ref{r-act-2}) is equivalent to
\be
\gamma\cdot_\chi [\rho]=[\rho]\;\; \mbox{for every}\; \gamma\in
\Gamma_\chi,
\label{ext-2-chi}
\ee
and hence if (\ref{r-act-2}) is satisfied,
since $\Gamma_\chi$ is cyclic,
we can apply Proposition \ref{cyclic} to conclude that
a semisimple representation of $K$ extends to a semisimple
representation of $G_\chi$.
To see this, let   $\Int^\chi_g$ the inner automorphism of $G_\chi$
defined by $g\in G_\chi$. 
We have to take some care about what we mean by the 
inner automorphism defined by $g$ since the underlying sets of  $G$
and $G_\chi$ coincide, but the group structures are different.
The restriction of $\Int^\chi_g$ to $K$ is  the homomorphism 
$$
\Int^\chi_g(h)=g\twist h\twist g^{-1}_\chi, \;\;\mbox{for}\;\; g\in
G_\chi \;\;\mbox{and}\;\;h\in K,
$$
where $g^{-1}_\chi$ denotes the inverse of $g$ with respect to the
operation $\twist$ defined by (\ref{twist-operation}).

By Proposition \ref{equivalence-action}, (\ref{ext-2-chi}) is equivalent to
\be
[\rho]=[\rho\circ \Int^\chi_g]\;\;\mbox{for every}\;\; g\in G_\chi.
\label{ext-1-chi}
\ee
Now, for every  $g\in G_\chi$ and $h\in K$
$$
\Int^\chi_g(h)=g\twist h\twist g^{-1}_\chi=
              g(h (g^{-1}_\chi)^{\chi^{-1}(h)})^{\chi^{-1}(g)}=
              gh^{\chi^{-1}(g)}g^{-1}
$$
since $\chi^{-1}(h)=1$ and $ (g^{-1}_\chi)^{\chi^{-1}(g)}=g^{-1}$.
We can thus conclude that 
$\Int^\chi_g=\Int_g\circ\chi^{-1}(g)$
and hence (\ref{ext-1-chi}) can be rewritten as 
$[\rho]=[\rho\circ\Int_g\circ\chi^{-1}(g)]$  
or equivalently
$$
[\rho\circ\chi(g)]=[\rho\circ\Int_g].
$$
Hence
$$
\gamma\cdot_\chi [\rho]=[\rho\circ \Int^\chi_{g_\gamma}]=
[\rho\circ \Int_{g_\gamma}\circ \gamma^{-1}],
$$
and (\ref{ext-2-chi}) is thus equivalent to (\ref{r-act-2}).

We have thus completed the proof of Theorem \ref{main-theorem}.
\qed

\remark 
Another way to prove our main theorem could be, perhaps,  
the following.
The twisted Higgs bundle $(E,\theta)$  defines
a Higgs bundle $(V,\Theta)$ over $X$  by taking
$V=\bigoplus_{i=0}^{n-1} E\otimes L^i$, where $n$ is the order of $\Gamma$,
and 
$$
\Theta=\left(\begin{array}{cccc}
0&...& 0& \theta\\
\theta &0& ...& 0\\
...&...& ... & ...\\
0&...&\theta & 0
\end{array}
\right)
$$
If $(E,\theta)$ is stable it should not be difficult to prove 
that $(V,\Theta)$ is polystable, defining then a representation
of $G$. This is the representation induced by the 
representation of $K$ corresponding to the twisted Higgs bundle $(E,\theta)$.
One would need then to characterize these representations of $G$  and
show that they are in bijection with the representations of $G_\chi$. 

\noindent {\bf Acknowledgements.}
The research carried out in this paper has taken place 
during visits of the second author to Universidad Aut\'onoma
(Madrid) and of the first author to TIFR  (Mumbai) and Ecole Polytechnique (Palaiseau),
 as well as joint visits
to ICTP (Trieste) and the Mathematical Institute (Oxford). We wish to thank these
institutions for their hospitality and support.
This work has been partly supported by the Spanish MEC PB95-0185. 
Both authors have benefited from workshops organised by the VBAC
research group under the European algebraic networks AGE and 
Europroj, supported by the European Union.

%%%%%%%%%%%%%%%%%%%%%%%%%%%%%%%%%%%%%%%%%%%%%%%%%%%%%%%%%%%%%%%%%%%%

%%%%%%%%%%%%%%%%%%%%%%

\noindent Departamento de Matem\'aticas, Universidad Aut\'onoma de Madrid,
28049 Madrid, Spain. \\
 {\tt oscar.garcia-prada@uam.es}\\

\noindent School of Mathematics, Tata Institute of Fundamental
Research, Homi Bhabha Road, Mumbai -- 400 005, India.\\
 {\tt ramanan@math.tifr.res.in}

%%%%%%%%%%%%%%%%%
\end{document}